\documentclass[12pt]{amsart}
\usepackage{amscd}
\newtheorem{theorem}{Theorem}
\newtheorem{lemma}{Lemma}
\newtheorem{proposition}{Proposition}
\newtheorem{corollary}{Corollary}
\theoremstyle{remark}

\newtheorem{definition}{Definition}
\newcommand{\T}{\mathbb T}
\newcommand{\Z}{\mathbb Z}
\newcommand{\R}{\mathbb R}
\newcommand{\mS}{\mathbb S}
\newcommand{\cH}{\mathcal H}
\newcommand{\cD}{\mathcal D}
\newcommand{\cT}{\mathcal T}
\newcommand{\cG}{\mathcal G}
\newcommand{\cU}{\mathcal U}
\newcommand{\cL}{\mathcal L}
\newcommand{\M}{\mathcal M}
\newcommand{\bq}{\mathbf{q}}
\newcommand{\bm}{\mathbf{m}}

\newcommand{\fui}{\varphi}
\newcommand{\ep}{\varepsilon}
\newcommand{\ga}{\gamma}
\newcommand{\Ga}{\Gamma}
\newcommand{\de}{\delta}
\newcommand{\lam}{\lambda}
\newcommand{\te}{\theta}
\newcommand{\lV}{\left\Vert}
\newcommand{\rV}{\right\Vert}
\newcommand{\rip}{\rangle}
\newcommand{\lip}{\langle}
\newcommand{\graph}{\hbox{graph\,}}
\begin{document}
\title[A minimax selector]{A minimax selector for a class of Hamiltonians on
cotangent bundles} 
\author[R. Iturriaga]{Renato Iturriaga}
\thanks{Partially supported by  CONACYT-M\'exico grant $\#$ 28489-E}
\address{CIMAT  \\
         A.P. 402, 36000 \\
         Guanajuato. Gto. \\
         M\'exico.}
\email{renato@cimat.mx}

\author[H. S\'anchez M.]{H\'ector S\'anchez-Morgado}
\address{Instituto de Matem\'aticas. UNAM  \\
         Ciudad Universitaria C. P. 04510\\
         M\'exico, DF\\
         M\'exico.}
\email{hector@matem.unam.mx}
\date{July 2000} 
\thanks{2000 Mathematics Subject Classification 37J10, 37J50} 

\begin{abstract}
We construct a minimax selector for eventually quadratic hamiltonians
on cotangent bundles. We use it to give a relation between Hofer's
energy and Mather's action minimizing function. We also study the
local flatness of the set of twist maps.
\end{abstract}
\maketitle
\section{Introduction}
One  purpose of this paper is to have some understanding of the
relation between two approaches to the study of hamiltonian systems,
namely the dynamical point of view of Mather \cite{M} and the
geometric point of view of Hofer and Zehnder \cite{hofer}.

The main result concerning this purpose is an extension to arbitrary 
cotangent bundles of a result
given by Siburg \cite{siburg2} for the cotangent bundle of an $n$
- torus. See Theorem \ref{mather} below.

To do so we follow the suggestion of   Bialy and Polterovich
\cite{B-P} and construct a  minimax selector for a certain class 
of Hamiltonians.
They conjectured that a selector  could be constructed 
for symplectic manifolds admitting a nice Floer homology. This program
has been carried out by Schwarz \cite{schwarz} for symplectic compact
manifolds using Floer homology. 

Since our interest is mainly on convex superlinear Hamiltonians
defined on cotangent bundles, we found very appealing to use the
methods developed by  Gol\'e \cite{gole} to find periodic orbits 
for Hamiltonians that are not necessarily convex but quadratic outside
of a neighbourhood of the zero section.

More precisely, let  $(M,\lip,\rip)$  be a compact connected
orientable riemannian manifold, such that there are no nontrivial
contractible closed geodesics. 
Let $\pi:T^*M\to M$ be its cotangent bundle endowed 
with its natural symplectic structure $\Omega=-d\lam$,
and let  
\[U_r=\{(q,p):|p|\le r\}.\]

Let ${\cH}_0$ be the set of smooth functions on
$T^*M\times{\mS}^1$ with support contained in $U_1\times{\mS}^1$. 
Let $\cH_R$ be the set of smooth functions on
$T^*M\times{\mS}^1$ such that $H=\frac R2 (|p|^2-1)$ for $|p|\ge 1$. 
For $H\in {\cH}_0\cup {\cH}_R$, let $\fui_t(x)=\fui^H_t(x)$ be
the solution to the Hamiltonian system $\dot{x}= X_H(x,t)$ with 
$\fui_0(x)=x$ and define $\fui_s^t(x)=\fui_t\circ\fui_s^{-1}(x)$.

Letting
\[{\cD}_0=\{\fui^H_1:H\in{\cH}_0\},\quad{\cD}_R=\{\fui^H_1:H\in{\cH}_R\},\]
for any $\psi\in\cD_R$ we have
${\cD}_R=\{\psi\circ\phi:\phi\in{\cD}_0\}$.

For $H\in {\cH}_0\cup {\cH}_R$, let 
\[E^+(H)=-\int_0^1\min_{x\in U_1}H(x,t)dt,\quad
E^-(H)=-\int_0^1\max_{x\in U_1}H(x,t)dt,\] and 
\[\lV H_t\rV=\max_{x\in U_1}H(x,t) -\min_{x\in U_1}H(x,t).\]

If $H\in {\cH}_0\cup\cH_R$, then $H=0$ on $\partial U_1$ and so 
$E^+(H)\ge 0\ge E^-(H)$.

For $H\in\cH_R$, define the {\em contractible action spectrum}
\[\sigma_c(H)=\left\{\int_{\Ga}\lambda -Hdt:\Ga\,
\mbox{is a contractible 1-periodic orbit of}\, \fui_t^H\right\}\]

\begin{theorem}\label{critical}
For any nonzero $v\in H^*(M,\R)$ and $H\in\cH_R$ there is
$c_v(H)\in\R$ such that 
\begin{itemize}
\item[1.] $c_v(H)\in\sigma_c(H)$
\item[2.] For  $K\in\cH_R$,  $E^-(K-H)\le c_v(K)-c_v(H)\le E^+(K-H)$.
\end{itemize}
\end{theorem}

In particular, for $H\in\cH_R$ define the selectors
\[\ga_-(H) = c_1(H),\, \ga_+(H) =c_{[\mu]}(H)\]
where $\mu$ is the orientation form.
These selectors have the following key properties:

\begin{theorem}\label{selector} Let $H,K\in\cH_R$. Then
\begin{itemize}
\item[1.] If $\fui_1^H=\fui_1^K$ then $\ga_\pm(H)=\ga_\pm(K)$.
\item[2.] If $H< 0$ in the interior of $U_1$ then $\ga_-(H)> 0$.
\item[3.] $E^-(H) \le\ga_-(H)\le\ga_+(H)\le E^+(H)$.
\item[4.] If $H\le 0$ on $U_1$ then $\ga_-(H)\ge 0$.
\end{itemize}
\end{theorem}

If  $\fui$ is generated by $H\in\cH_R$ i.e. $\fui^H_1=\fui$,
define $\ga_\pm(\fui)=\ga_\pm(H)$. 

For $\fui\in\cD_R$ define its energy as 
\[E(\fui)=\inf\left\{\int_0^1\lV H_t\rV
  dt:H\in{\cH}_R,\fui^H_1=\fui\right\}.\]

For $\phi\in{\cD}_0$  let 
\[\lV \phi \rV=\inf\{\int_0^1\lV K_t\rV dt:K\in{\cH}_0,
\fui^K_1=\phi\}\]

If $\fui,\psi\in\cD_R$ define $d(\fui,\psi)=\lV\psi^{-1}\fui\rV$.

As proved by Lalonde and Mc Duff in \cite{L-MD}, $\lV\fui \rV=0$ only
when $\fui$ is the identity and then $d$ is  a metric in $\cD_R$.

For $L:TM\times{\mS}^1\to{\R}$ a convex superlinear Lagrangian with
complete Euler-Lagrange flow, let
$\beta=\beta_L:H_{1}(M,\R)\rightarrow \R$ be the Mather's beta function.
The following theorem gives a connection between  the
points of view of Mather and Hofer-Zehnder.

\begin{theorem}\label{mather}
If $\fui$ is generated by a convex Hamiltonian $H \in \cH_R$, then 
  \[\beta_L(0)\le E(\fui) \]
where $L$ is the Legendre transform of $H$.
\end{theorem}

A {\it twist map} $F$ is a diffeomorphism of a neighborhood
$U$ of the zero section in $T^*M$ onto itself satisfying the following:
\begin{itemize}
\item If $F(q,p)=(Q,P)$, then the map $\Phi(q,p)=(q,Q)$ is an 
embedding of $U$ in $M\times M$
\item  F is {\it exact symplectic}, that is, $F^*\lam-\lam = 
d (S\circ\Phi)$ 
for some  real function $S$ on $\Phi(U)$.
\end{itemize}
The function $S$ is called a {\it generating function} for $F$.

Fix $r>0$ smaller than the injectivity radius of the given metric. 
Let $\cT$ be the set of twist maps $\fui|U_1$ with $\fui\in\cD_r$.
Any map in $\cT$ has a generating function satisfying 
\begin{itemize}
\item $\dfrac{\partial^2 S}{\partial q\partial Q}$ is negative definite.
\item $S(q,Q)=\dfrac{d(q,Q)^2}{2r}$ for $d(q,Q)\ge r$.
\end{itemize}

The following proposition gives the local flatness of Hofer's metric
for twist maps and is a generalization of Theorem 1 in Siburg's paper
\cite{siburg1}.

\begin{proposition}\label{flat}
For any $\fui\in{\cT}$ there is a $C^1$ neighborhood $\mathcal O$
of $\fui$ such that if $\fui_0,\fui_1\in\mathcal O$ and $S_0,S_1$
are their generating functions as above, then
\[d(\fui_0,\fui_1)=\lV S_0-S_1\rV.\]
\end{proposition}

In the appendix  we consider a general convex superlinear 
     Hamiltonian $H:T^*M\times{\mS}^1\to\R$ and its Legendre transform
     $L:TM\times{\mS}^1\to{\R}$.
     Let  $\alpha:H^{1}(M,\R)\rightarrow \R$ be the convex dual
     of $\beta$.

     Recall that the image of a one-form in $M$ is a lagrangian submanifold
     of $T^*M$ if and only if the form is closed. 
     Following Bialy and Polterovich \cite{B-P} we define

     \[c(H)=\inf\left\{\int_0^1\max_{q\in M}H(q,\te(q),t)dt: \te 
     \text{ is a closed 1-form }\right\}.\]
 
   More generally we can define for any cohomology class 
   $[\omega]\in H^1(M,\R)$
   \[c(H,[\omega])=\inf\left\{\int_0^1\max_{q\in M}H(q,\te(q),t)dt: 
     \te\in[\omega] \right\}.\]

   When $H$ is autonomous this reduces to 
   \[c(H,[\omega])=\inf_{\te\in[\omega]}\max_{q\in M}H(q,\te(q))\] 
   and it was proved in \cite{gafa} that
   $c(H,[\omega])=\alpha([\omega])$. 
   In the non-autonomous case we recover the following inequality

   \begin{proposition}
     \label{lagrange}
     \begin{equation}
       \label{eq:desgafa}
       \alpha([\omega]) \le c(H,[\omega]).
     \end{equation}
   \end{proposition}

   Taking  infimum on both sides of inequality \eqref{eq:desgafa},
   we obtain 

   \begin{corollary}
     \label{gafa}
     \[-\beta(0)\le c(H).\]
   \end{corollary}

   For the last application we go back to our special family $\cH$
   \begin{corollary}\label{equal}
     Let $H\in{\cH}$ be convex and suppose there is a lagrangian 
     section in $U_1$
     consisting of fixed points of $\fui^H_t$. Then
     \[E(\fui^H_1)=\beta (0)=-c(H).\]
   \end{corollary}

 The plan of the article is as follows, in section \ref{generating}
    we introduce the generating families whose critical values belong
    to the contractible action spectrum. In section \ref{values} we
    prove theorems \ref{critical} and \ref{selector}. In section
    \ref{main} we recall the definition of Mather's beta function and
    prove theorem \ref{mather}.   In section \ref{hofer} we apply the
    results of Lalonde and Mc Duff \cite{L-M} describing the geodesics
    of Hofer's metric for certain classes of symplectic manifolds that
    include cotangent bundles of compact manifolds and prove
    proposition \ref{flat}.

    We thank K.F. Siburg for his useful comments on a first draft of
    this paper.

\section{Generating families}
\label{generating}

Let  $g_t$ be the geodesic flow on $T^{\ast}M$,  $d$ be the distance in $M$ 
defined by the metric and fix $r>0$ smaller than the injectivity radius.
The map $g_r|U_1$ is twist since $G_r(q,p)=(q,\exp_q(r p))$
is an embedding of $U_1$ in $M\times M$ and $S(q,Q) = {d(q,Q)^2}/{2r}$
is a generating function. Write $g_{-r}=(Q_g,P_g)$.

For $H\in{\cH}_R$, let $\fui_t(x)=\fui^H_t(x)$.
Since the $C^1$ distance from $\fui_s^t$ to the identity is $O(s-t)$ on $U_1$,
there is $N$ such that 
\[\Phi_s^t(q,p)=(q,Q_g\circ\fui_s^t(q,p))\]
is an embedding of $U_1$ in $M\times M$ for $|t-s|\le 1/N$.
For $0\le j< N$, let $s_j=j/N$,
$$\fui^j=\fui_{s_j}^{s_{j+1}}, \, F_{2j}=g_{-r}\circ\fui^j,\, 
F_{2j+1}=g_r, \,\Phi_{2j}=\Phi_{s_j}^{s_{j+1}}, \, \Phi_{2j+1}=G_r$$
Let $\Ga_x^j(t)=\fui_{s_j}^t(x),t\in[s_j,s_{j+1}]$, and
define the {\em action} $A_j:U_1\to{\R}$ by
\begin{equation}\label{action}
A_j(x)=\int_{\Ga_x^j}\lam-H dt,
\end{equation}
then $dA_j=(\fui^j)^\ast\lambda-\lambda$. Thus, defining
$$S_{2j}=
(A_j-S\circ G_r\circ F_{2j})\circ\Phi_{2j}^{-1},$$
we have
\begin{equation}\label{eq:generating}
\begin{split}
d\,S_{2j}\circ \Phi_{2j} & 
=dA_j - F_{2j}^\ast\,d\,S\circ G_r\\
& =(\fui^j)^\ast\lambda-\lambda-(\fui^j)^\ast g_{-r}^\ast
(g_r^\ast\lambda-\lambda)\\
& = F_{2j}^\ast\lambda-\lambda.
\end{split}
\end{equation}
 
Therefore $F_{2j}$ is a twist map with generating function $S_{2j}$. 
Let $S_{2j+1}=S$ and  
\[{\M}_N=\{(q_0,q_1,\ldots,q_{2N}):(q_j,q_{j+1})\in\Phi_j(U_1), q_{2N}=q_0\}.\] 
The function $W_N:{\M}_N\to{\R} $ given by
\[W_N:{\bq}=(q_0,q_1,\ldots,q_{2N})
\mapsto\sum_{j=0}^{2N-1}S_j(q_j,q_{j+1}).\]
is called a generating family.
Letting $(q_j,p_j)=\Phi_j^{-1}(q_j,q_{j+1})$, $p_{2N}=p_0$,
$(q_{j+1},P_{j+1})=F_j(q_j,p_j)$, $P_{2N}=P_0$, we have by \eqref{eq:generating}
\begin{equation}\label{gradient}
dW_N=\sum_{j=0}^{2N-1}(P_j-p_j)dq_j.
\end{equation}

Therefore ${\bq}$ is a critical point of $W_N$ if and only if 
$(q_{j+1}, p_{j+1})=F_j(q_j,p_j)$ for $0\le j<2N$, so 
$(q_{2j},p_{2j})=\fui_{j/N}=(q_0,p_0)$, $(q_0,p_0)$
is a fixed point of $\fui_1$, and $\Ga(t)=\fui_t(q_0,p_0)$ is a 1-periodic 
orbit in $U_1$. In such a case
\[W_N({\bq})=\sum_{j=0}^{N-1}A_j(q_j,p_j)=\int_{\Ga}\lam-Hdt.\]

\begin{proposition}\label{derivatives}
Given $\ep>0$ there is $N(\ep)$ such that for $N\ge N(\ep)$
we have
\begin{equation}
  \label{eq:odd-close}
\lV  \frac{\partial P_{2j+1}}{\partial q_{2j+1}}-
\frac{\partial p_{2j+1}}{\partial q_{2j+1}} \rV <\ep
\end{equation}
for any $j$ such that
$F_{2j}\circ\Phi_{2j}^{-1}(q_{2j},q_{2j+1})=G_r^{-1}(q_{2j+1},q_{2j+2})$,
and
\begin{equation}
  \label{eq:even-close}
\lV  \frac{\partial P_{2j}}{\partial q_{2j}}-
\frac{\partial p_{2j}}{\partial q_{2j}} \rV <\ep
\end{equation}
for any $j$ such that 
$g_r\circ G_r^{-1}(q_{2j-1},q_{2j})=\Phi_{2j}^{-1}(q_{2j},q_{2j+1})$.
\end{proposition}
{\it Proof}
\begin{itemize}
\item [A.] From 
$(q_{2j+1},P_{2j+1})=F_{2j}\circ\Phi_{2j}^{-1}(q_{2j},q_{2j+1})$,
$F_{2j}=g_{-r}\circ\fui^j$ one gets
\[  \frac{\partial P_{2j+1}}{\partial q_{2j+1}}=
\frac{\partial (P_g\circ\fui^j)}{\partial p}\circ\Phi_{2j}^{-1}
\left(\frac{\partial (Q_g\circ\fui^j)}{\partial p}\circ\Phi_{2j}^{-1}
\right)^{-1}. \]

\item [B.] From $(q_{2j+1},p_{2j+1})=G_r^{-1}(q_{2j+1},q_{2j+2})$ one gets
\[  \frac{\partial p_{2j+1}}{\partial q_{2j+1}}=
\frac{\partial P_g}{\partial p}\circ g_r\circ G_r^{-1}
\left(\frac{\partial Q_g}{\partial p}\circ g_r\circ G_r^{-1}
\right)^{-1}. \]

\item [C.] From $q_{2j+1}=Q_g\circ\fui^j(q_{2j},p_{2j})$ one gets
\[ \frac{\partial p_{2j}}{\partial q_{2j}}=
-\left(\frac{\partial(Q_g\circ\fui^j)}{\partial p}\circ\Phi_{2j}^{-1}
\right)^{-1}\frac{\partial(Q_g\circ\fui^j)}{\partial q}\circ\Phi_{2j}^{-1}.\]

\item [D.] From $q_{2j-1}=Q_g(q_{2j},P_{2j})$ one gets

\[  \frac{\partial P_{2j}}{\partial q_{2j}}=
- \left(\frac{\partial Q_g}{\partial p}\circ g_r\circ G_r^{-1}
\right)^{-1}\frac{\partial Q_g}{\partial q}\circ g_r\circ G_r^{-1}.\]
\end{itemize}

Suppose that 
$F_{2j}\circ\Phi_{2j}^{-1}(q_{2j},q_{2j+1})=G_r^{-1}(q_{2j+1},q_{2j+2})$
so that  
\[\fui^j\circ\Phi_{2j}^{-1}(q_{2j},q_{2j+1})=
g_r\circ G_r^{-1}(q_{2j+1},q_{2j+2}),\]
then
\begin{multline*}
D(g_{-r}\circ\fui^j)(\Phi_{2j}^{-1}(q_{2j},q_{2j+1}))=\\
Dg_{-r}(g_r\circ G_r^{-1}(q_{2j+1},q_{2j+2})) 
D\fui^j(\Phi_{2j}^{-1}(q_{2j},q_{2j+1}))
\end{multline*}
Suppose
$g_r\circ G_r^{-1}(q_{2j-1},q_{2j})=\Phi_{2j}^{-1}(q_{2j},q_{2j+1}),$
then
\[D(g_{-r}\circ\fui^j)(\Phi_{2j}^{-1}(q_{2j},q_{2j+1}))=
D(g_{-r}\circ\fui^j)(g_r\circ G_r^{-1}(q_{2j-1},q_{2j})).\]

Taking $N$ sufficiently large, one makes any $\fui^j|U_1$ as $C^1$- close 
to the identity as one wants. Thus, \eqref{eq:odd-close} 
follows from items A, B, and \eqref{eq:even-close} follows from items C,D.

\qed

As observed by Gol\'e \cite{gole} one
can define a path $\sigma_k(q,Q)$ between $q$ and 
$Q\in\Phi_k(U_1\cap T^*_qM)$ that for $k$ odd concides with the unique
geodesic between these points.  Therefore to each ${\bq}\in{\M}_N$ 
we associate a polygonal loop $c(\bq)$.  We can work with the
component ${\M}_N^e$ of ${\M}_N$ consisting of points $\bq$ such that
$c({\bq})$ is homotopically trivial. Thus, $\sigma_c(H)$ is the set of
critical values of $W_N:{\M}_N^e\to \R$ and so it is nowhere dense.

Gol\'e  proved that the gradient flow of $W_N$ has 
the index pair $(B_N,B_N^-)$,
\[B_N=\{\bq\in{\M}_N^e:d(q_j,q_{j+1})\le |a_j|\},\]
\[B_N^-=\{\bq\in B_N:d(q_j,q_{j+1})= |a_j| \quad
\text{for some even}\quad j \},\]
where 
\[a_j=\begin{cases} r & j \,\text{odd} \\ 
\frac RN -r & j \,\text{even} \end{cases}\]

\begin{lemma}\label{floer} {\em Floer} \cite{floer}. 

Let $\zeta^\tau$ be a one parameter family of flows on a manifold $\M$.
Suppose that $\Sigma^0$ is a compact submanifold invariant under $\zeta^0$.
Assume moreover that $\Sigma^0$ is {\em normally hyperbolic} i. e. there is
a decomposition
\[T_{\Sigma^0}\M = T\Sigma^0\oplus E^+\oplus E^-\]
invariant under the covariant linearization of the vector field $V_0$
generating $\zeta^0$ with respect to some metric  $\lip,\rip$, so that for 
some $m>0$:
\[\lip\xi,DV_0\xi\rip\le -m\lip\xi,\xi\rip,\quad \xi\in E^-\]
\[\lip\xi,DV_0\xi\rip\ge m\lip\xi,\xi\rip,\quad \xi\in E^+\]
Suppose that there is a retraction $\alpha:\M\to\Sigma^0$ and that there is
an index pair $(B,B^-)$ for all $\zeta^\tau$. Then 
\begin{itemize}
\item 
There is $u\in H^*(B,B^-)$ of dimension $\dim E^+$ so that the map
\[T:H^*(\Sigma^0)\to H^*(B,B^-),\,T(v)= (\alpha|_B)^*v\cup u\]
is an isomorphism. 
\item
If $\Sigma^\tau$ denotes the maximal invariant set for the flow $\zeta^\tau$, 
the homomorphism in \v{C}ech cohomology
$(\alpha|\Sigma^\tau)^*:H^*(\Sigma^0)\to H^*(\Sigma^\tau)$ is injective.
\end{itemize}
\end{lemma}

\section{Critical values of the action}\label{values}

 Let $H\in\cH_R$, and define  $H_0=\frac R2 (|p|^2-1)$.  Let
    $H_\tau=H_0+\tau(H-H_0)$ and $\fui^\tau_t$ be the corresponding
    flow. Let $N$ be sufficiently large to have a  decomposition of
    all $\fui^\tau_1$ in $2N$ twist maps. Let $W_N^\tau$ be the
    corresponding generating family and $\zeta^\tau$ its gradient flow.

    \begin{lemma}{\em Gol\'e} \cite{gole}.  Let
    $\Sigma^0=\{\bq\in{\M}_N:q_k=q_0,\,\forall k\}$.  Then $\Sigma^0$
    is a normally hyperbolic invariant set for $\zeta^0$ and it is a
    retract of ${\M}_N^e$  \end{lemma}

    By Lemma \ref{floer},  $T_N:H^*(M)\to H^{*+k}(B_N,B_N^-)$ is an
    isomorphism.

     Let  $B_N^a=\{\bq\in B_N:W_N(\bq)\le a\}$,
    $j_a:B_N^a\hookrightarrow B_N$ and consider the induced map
    $j_a^*:H^*(B_N,B_N^-)\to H^*(B_N^a,B_N^-)$.

    For $N$ of the form $N=2^m$, we follow Viterbo \cite {viterbo}
    and define  for any nonzero  $v\in H^*(M)=H^*(\Sigma^0)$
   
\[c_v(H,m)=\inf\{a:j_a^*T_N(v)\ne 0\}.\]
so that $c_v(H,m)$ is a critical point of $W_N$.

\begin{proposition}
There is $\bm$ such that for any nonzero $v\in H^*(M)$ and  $m\ge\bm$, 
$c_v(H,m)=c_v(H,\bm)$.
\end{proposition}

{\it{Proof}}.
Represent the points of ${\M}_{2N}$ in the form
\[\bq=(\eta_0,\xi_0,\xi_1,\eta_1,
\ldots,\eta_{2N-2},\xi_{2N-2},\xi_{2N-2},\eta_{2N-1},\eta_{2N}=\eta_0)\]
so that

\begin{multline*}
W_{2N}({\bq})=W_{2N}(\eta,\xi)= \\ \sum_{j=0}^{N-1}S_{4j}(\eta_{2j},\xi_{2j})
+S(\xi_{2j},\xi_{2j+1})+S_{4j+2}(\xi_{2j+1},\eta_{2j+1})
+S(\eta_{2j+1},\eta_{2j+2}).
\end{multline*}
Let $h_j=
G_r\circ g_{-R}\circ\fui_{j/N}^{2j+1/2N}\circ(\Phi_{j/N}^{j+1/N})^{-1}$
and define $h:{\M}_N\to M^{2N}$ by
$h(\eta)=(h_0(\eta_0,\eta_1),\ldots,h_{N-1}(\eta_{2N-2},\eta_{2N-1})).$
Then $W_{2N}(\eta,h(\eta))=W_N(\eta)$ and
\[\frac{\partial W_{2N}}{\partial\xi_j}(\eta,h(\eta))=0.\]

Define $f:{\M}_{2N}\to\R$ by $f(\eta,\xi)=W_{2N}(\eta,\xi)-W_N(\eta)$,
then
\begin{eqnarray}
  \label{eq:deriva-W}
    \frac{\partial f}{\partial\xi_{2j+i}} & = &
\frac{\partial W_{2N}}{\partial\xi_{2j+i}}=
P_{4j+1+i}-p_{4j+1+i},\,i=0,1\\
\frac{\partial f}{\partial\eta} & = &
\frac{\partial W_{2N}}{\partial\eta}-
\frac{\partial W_{2N}}{\partial\eta}(\eta,h(\eta))-
\frac{\partial W_{2N}}{\partial\xi}Dh.
\end{eqnarray}
Thus $Df|\graph h=0$ and

\[\dfrac{\partial^2 f}{\partial(\xi_{2j},\xi_{2j+1})}
=\begin{bmatrix}
\dfrac{\partial(P_{4j+1}-p_{4j+1})}{\partial\xi_{2j}} &
\dfrac{\partial^2 S}{\partial q\partial Q}(\xi_{2j},\xi_{2j+1})\\
\dfrac{\partial^2 S}{\partial q\partial Q}(\xi_{2j},\xi_{2j+1}) &
\dfrac{\partial(P_{4j+2}-p_{4j+2})}{\partial\xi_{2j+1}}
\end{bmatrix}.\]

Since $\dfrac{\partial^2 S}{\partial q\partial Q}$ is negative
definite on $G_r(U_1)$, we have that 

$\dfrac{\partial^2 f}{\partial(\xi_{2j},\xi_{2j+1})}|\graph h$ 
is invertible if
$\dfrac{\partial(P_{4j+1+i}-p_{4j+1+i})}{\partial\xi_{2j+i}}|\graph h$,
$i=0,1$  are sufficiently small.

We know from Proposition \ref{derivatives} that this holds if 
$N=2^m$ is sufficiently large. In such a case, $\graph h$ is a
nondegenerate critical manifold of $f$. By the generalized Morse's
lemma there is a tubular neighborhood $\psi:E\to{\M}_{2N}$, 
with $E=E^+\oplus E^-$ a vector bundle over ${\M}_N$, such that
$f(\psi(\eta,\zeta))=|\zeta_-|^2-|\zeta_+|^2$
and so
\begin{equation}
  \label{eq:morse}
  W_{2N}(\psi(\eta,\zeta))=W_N(\eta)+|\zeta_-|^2-|\zeta_+|^2.
\end{equation}

Consider the commutative diagram
\[\begin{matrix}&T_{2N}\\&\nearrow\\&\\H^*(M)&\\&
\\&\searrow\\&T_N\end{matrix}
\begin{CD} 
H^{*+k+l}(B_{2N},B_{2N}^-) @>j_a^*>> H^{*+k+l}(B^a_{2N},B_{2N}^-)\\
@VV\psi^* V    @VV\psi^* V\\
H^{*+k+l}(\psi^{-1}(B_{2N},B_{2N}^-)) @>J_a^*>>  
H^{*+k+l}(\psi^{-1}(B^a_{2N},B_{2N}^-)) \\
@AA\text{Thom}A         @AA\cong A\\
H^{*+k}(B_N,B_N^-) @>j_a^*>> H^{*+k}(B^a_N,B_N^-)
\end{CD}  \]
where the upwards isomorphisms are as in \cite{viterbo}. 

The maximal invariant set $\Sigma_a$ for the gradient flow of $W_{2N}$
in $B^a_{2N}-B_{2N}^-$ consists of critical points and heteroclinic orbits. 
The critical points belong to $\graph h$ by \eqref{eq:deriva-W} and so
do the heteroclinic orbits by \eqref{eq:morse}. Thus
$(B^a_{2N},B_{2N}^-)$ and $(B^a_{2N},B_{2N}^-)\cap\psi(E)$ are index
pairs for $\Sigma_a$ and then
\[\psi^*:H^*(B^a_{2N},B_{2N}^-)\cong 
H^*(\psi^{-1}(B^a_{2N},B_{2N}^-)).\] 

Therefore $j_a^*T_{2N}(v)=0$ if and only if $j_a^*T_N(v)=0$ and so
\[c_v(H,m) = c_v(H,m+1).\]

\qed

We consider in particular 
\[\ga_-(H) = c_1(H),\,\ga_+(H) = c_{[\mu]}(H).\]

\begin{lemma}\label{viterbo}{\em Viterbo} \cite{viterbo}.
Let $S_\tau$ be a smooth family of smooth functions.
Let $c(\tau)=S_\tau(x_\tau)$ be a critical value obtained
by minimax as above.
Assume that $dS_\tau(x)=0$ implies 
$\dfrac{\partial S_\tau(x)}{\partial \tau}\ge 0\: (\le 0)$. 
Then $c(\tau)$ is increasing (decreasing).
\end{lemma}

The following Lemma completes the proof of Theorem \ref{critical}
\begin{lemma}\label{fundamental} 
For $v\in H^*(M)$, $H,K\in\cH_R$
\[E^-(K-H)\le c_v(K)-c_v(H)\le E^+(K-H).\]
\end{lemma}

{\it{Proof}}. For $\tau\in [0,1]$, define $H_\tau= H+ \tau (K-H)$, 

Let $S_\tau=W_\tau-E^-(K-H)\tau$ then $c_v(H_\tau)-E^-(K-H)\tau$ 
is a critical value of $S_\tau$ and
\[c_v(H_\tau)= W_\tau(x_\tau)=\int_{S^1}\psi_\tau^*(\lam-H_\tau dt)\]
where $\psi_\tau(t)=\fui_t^\tau(x_\tau)$.
Then  
\begin{equation*}
\begin{split}
  \frac{\partial W_\tau}{\partial\tau} (x_\tau)& = 
\int_{S^1}d\psi_\tau^*\lam(\frac{\partial\psi_\tau}{\partial\tau}) + 
\psi_\tau^*i(\frac{\partial\psi_\tau}{\partial\tau})d\lam\\
&-\int_0^1\left(\frac{\partial H_\tau}{\partial\tau}(\psi_\tau(t),t)+
dH_\tau(\psi_\tau(t),t)\frac{\partial\psi_\tau}{\partial\tau}\right)dt\\
&  =  \int_0^1(H-K)(\fui_t^\tau (x_\tau ),t)dt\ge E^-(K-H).
\end{split}
\end{equation*}
By Lemma, \ref{viterbo} $c_v(H_\tau)-E^-(K-H)\tau$ is increasing.
Similarly $c_v(H_\tau)-E^+(K-H)\tau$ is decreasing. Then
\[c_v(K)-E^+(K-H)\le c_v(H)\le c_v(K)-E^-(K-H)\]
\qed

\begin{corollary}\label{monotone}
If $H\le K$  then $c_v(H)\ge c_v(K)$ for any nonzero $v\in H^*(M)$.
\end{corollary}

{\it Proof of Theorem} \ref{selector}.

\begin{lemma}\label{basic}
If $H,K\in\cH_R$ are such that the corresponding time one maps $\fui_1$ and  
$\psi_1$ are equal then $\sigma_c(H)=\sigma_c(K)$.
\end{lemma}

{\it Proof}. 
Let $x_0$ in the border of $U_1$.
Let $x_1$ in $\mbox{Fix}(\fui_1)=\mbox{Fix}(\psi_1)$. 
Then $\fui_t(x_0)=\psi_t(x_0)$ for
all $t$. Let $\beta $ be a curve on $T^*M$ such that $\beta  (0)=x_0$
and $\beta  (1)=x_1$. Define $\sigma_1,\sigma_2:[0,1]^2\to T^*M\times\R$ by
\begin{eqnarray*}
  \sigma_1(s,t) & = & (\fui_t(\beta (s)),t)\\
  \sigma_2(s,t) & = & (\psi_t(\beta (s)),t).
\end{eqnarray*}

Then

$$0=\int_{\sigma_1}d(\lambda -Hdt)=
\int_{\fui_t(x_1)}\lambda -Hdt- \int_{\fui_t(x_0)}\lambda -Hdt -
\int_{\fui_1(\beta )}\lambda+ \int_{\beta }   \lambda $$
 and 

$$0=\int_{\sigma_2}d(\lambda -Kdt)=
\int_{\psi_t(x_1)}\lambda -Kdt- \int_{\psi_t(x_0)}\lambda -Kdt-
\int_{\psi_1(\beta )}\lambda+ \int_{\beta }   \lambda $$

So
$$\int_{\fui_t(x_1)}\lambda -Hdt=\int_{\psi_t(x_1)}\lambda -Kdt, $$

It remains to prove that $\fui_t(x_1)$ is contractible if and only if 
$\psi_t(x_1)$ is. To see that, define the path
 \[h_t(x)=\begin{cases} \fui_t(x) & \text{if}\quad t\in[0,1]\\
                        \psi_{2-t}(x) & \text{if}\quad t\in[1,2]
           \end{cases} .     \]
Let $\Omega(U_1)$ be the free loop space of $U_1$, and
define the continuous function $\chi:U_1\to \Omega(U_1)$ by 
$\chi(x)=(h_t(x))_{t\in[0,2]}$.
Since $M$ is connected all elements of $\chi(U_1)$ are homotopic to 
$\chi(x_0)$ which is homotopically trivial. Thus $\fui_t(x_1)$ is 
contractible if and only if $\psi_t(x_1)$ is and then
\[\sigma_c(H)=\sigma_c(K).\] 
\qed

Let $H\in C^\infty(T^*M\times{\mS}^1\times [0,1])$ be such that
$H_s(p,t)=H(p,t,s)$  belongs to $\cH_R$ for any $s\in[0,1]$. 
Assume now that $\fui^{H_s}_1=\psi$ for all $s\in[0,1]$, then
$\ga_\pm(H_s)\in\sigma_c(H_s)=\sigma_c(\psi)$.
Since $\sigma_c(\psi)$ is nowhere dense, the continuous function
$s\mapsto\gamma_\pm(H_s)$ must be constant and so
$\ga_\pm(H_0)=\ga_\pm(H_1)$. 

For $K\in\cH_R$, $L\in{\cH}_0$ let
$K\# L(p,t)=K(p,t)+L((\fui^K_t)^{-1}(p),t)$, then $H=K\# L\in\cH_R$ and
$\fui^H_t=\fui^K_t\fui^L_t$. Taking $H_0=\frac R2(|p|^2-1)$,
we can write any $H\in\cH_R$ as $H=H_0\# L$, where
$L(p,t)=(H-H_0)(g_t(p),t)$ belongs to ${\cH}_0$.

{\it Item} 1. We follow the argument of Hofer and Zehnder.

Let $H,K\in\cH_R$ such that $\fui^H_1=\fui^K_1$. Write 
$H=H_0\# L, K=H_0\# F$ with $L,F\in{\cH}_0$ and let
$\psi=\fui_1^L=\fui_1^F\in{\cD}_0$. 
By reparametrizing the time, one first homotope  $\fui^L_t$ to an arc $\psi_t$
in $\cD_0$ that is the identity for $t\in[0,1/4]$ and is equal to $\psi$
for $t\in[3/4,1]$, and do the same for $\fui^F_t$. Hence we can assume
that $L(p,t)=F(p,t)=0$ for $|t|_{\text{mod}\,1}<1/4$.
For $0<s\le 1$ one defines 
\[L_s(p,t)=s L(p/s,t),\, F_s(p,t)=s F(p/s,t).\]
Then
\[\fui^{L_s}_t(p)=s\fui^L_t(p/s),\, \fui^{F_s}_t(p)=s\fui^F_t(p/s).\]
For $t\in[3/4,1]$ one has 
\[\fui^{L_s}_t(p)=s\fui^L_t(p/s)=s\fui^F_t(p/s)=\fui^{F_s}_t(p).\]
Take $\beta:[3/4,1]\to[0,1]$ a smooth which is 0 near to 3/4 and is 1
near to 1. Define
\[\fui_{s,t}(p)=\begin{cases} \fui^{L_s}_t(p) & t\in[0,3/4] \\
(s+(1-s)\beta(t))\psi((s+(1-s)\beta(t))^{-1}p) & t\in[3/4,1] 
\end{cases} \]
and similarly $\psi_{s,t}$ by replacing $L_s$ for $F_s$. Then  
\[\fui_{s,1}=\psi_{s,1}=\psi\quad 0<s\le 1.\]
If $\hat{L}_s, \hat{F}_s$ generate $\fui_{s,t},\psi_{s,t}$
respectively we have $\ga_\pm(H)=\ga_\pm(H_0\#\hat{L}_s)$, 
$\ga_\pm(K)=\ga_\pm(H_0\#\hat{F}_s)$. Note that 
\[\hat{L}_s(p,t)=\hat{F}_s(p,t)\quad t\in[3/4,1]\]
\[\hat{L}_s(p,t)=L_s(p,t)\quad t\in[0,3/4]\]
\[\hat{F}_s(p,t)=F_s(p,t)\quad t\in[0,3/4].\]
Thus
\begin{multline*}
|\ga_\pm(H)-\ga_\pm(K)|=|\ga_\pm(H_0\#\hat{L}_s)-\ga_\pm(H_0\#\hat{F}_s)|
\\\le E^+(\hat{L}_s-\hat{F}_s)- E^-(\hat{L}_s-\hat{F}_s)
\le s(E^+(L)-E^-(L)+E^+(F)-E^-(F))
\end{multline*}
for all $0<s\le 1$ and then $\ga_\pm(H)=\ga_\pm(K)$.

For any $\ep>0$ let ${\mathcal K}_{\ep}$ be the set of functions $f(|p|)$ 
in $\cH_R$ such that $f$ is convex and $f(s)=-\ep$ for $s\le 1-2\ep/R$.

We claim that $\ga_\pm(K)=\ep$ for any $K\in{\mathcal K}_{\ep}$ . 
In fact, the only contractible periodic orbits are the constants and
these have action $\ep $. 

{\it Item} 2. Since $H$ is negative for $|p|\le 1$ and equals 
$\frac{R}{2}(|p|^2-1)$ for $|p|\ge 1$, we can find $\ep>0 $ small 
enough and $K\in{\mathcal K}_{\ep}$ such that $K\ge H$ and then by
Corollary \ref{monotone}
we have $\ga_-(H)\ge \ga_-(K)> 0$.

{\it Item} 3. 
By Lemma \ref{fundamental}, for $K_\ep\in{\mathcal K}_{\ep}$ we have
\[E^-(H-K_\ep)\le \ga_\pm(H)-\ga_\pm(K_\ep)\le E^+(H-K_\ep).\]
Now let $\ep\to 0$.

{\it Item} 4. Since $H\le 0$, $E^-(H)= 0$. By item 3, $\ga_\pm(H)\ge 0$.
\qed

The following Corollary is an immediate consequence of item 3 of
Theorem \ref{selector}.

\begin{corollary}\label{corollary}
If $\fui\in{\cD}_R$
\[\ga_-(\fui)\le\ga_+(\fui)\le E(\fui),\quad 
\ga_+(\fui)-\ga_-(\fui)\le E(\fui).\]
\end{corollary}

\section{Proof of Theorem \ref{mather}}\label{main}

Let us recall the main concepts introduced by Mather in \cite{M}.
Let $L:TM\times{\mS}^1\to{\R}$ be a  convex superlinear Lagrangian
with complete Euler-Lagrange flow.
Let ${\M}(L)$ be the set of probabilities on the Borel
$\sigma$-algebra  of $TM$ that have compact support and are 
invariant under the Euler-Lagrange flow $\phi_{t}$.
Let $H_{1}(M,\R)$ be the first real homology group of $M$.
Given a closed one-form $\omega$ on $M$ and $\rho\in H_{1}(M,\R)$, 
let $<\omega,\rho>$ denote the integral of $\omega$ on any closed
curve in the homology class $\rho$.
If $\mu\in {\M}(L)$, its {\it rotation vector} is defined as the
unique $\rho(\mu)\in H_{1}(M,\R)$ such that
\[<\omega,\rho(\mu)>=\int\omega\,d\mu,\]
for all closed one-forms on $M$.
The integral on the right-hand side is with respect to $\mu$ with 
$\omega$ considered as a function $\omega:TM\rightarrow \R$.
The function $\rho:{\M}(L)\rightarrow H_{1}(M,\R)$ is
surjective \cite{M}. 
The {\it action} of $\mu\in {\M}(L)$ is defined by
\[A_{L}(\mu)=\int L\,d\mu.\]
Finally we define the function $\beta:H_{1}(M,\R)\rightarrow \R$ by
\[\beta(\gamma)=\inf\{A_{L}(\mu):\;\rho(\mu)=\gamma\}.\]
The function $\beta$ is {\it convex} and {\it superlinear} and the
infimum can be shown to be a {\it minimum} \cite{M} and the measures 
at which the minimum is attained are called {\it minimizing measures}.
In other words, $\mu\in {\M}(L)$ is a minimizing measure iff
\[\beta(\rho(\mu))=A_{L}(\mu).\]

Each contractible 1-periodic orbit is the support of an invariant
probability measure with zero homology and the same action. So
$$\beta(0)\le\ga_- (H).$$

By Corollary \ref{corollary}, $\ga_-(\fui)\le E(\fui)$,
and this concludes the proof of Theorem \ref{mather}.
\qed

\section{A Hofer's distance and twist maps}\label{hofer}
For $\phi\in{\cD}_0$  let 
\[\lV \phi \rV=\inf\{\int_0^1\lV L_t\rV dt:L\in{\cH}_0,
\fui^L_1=\phi\}.\]

If $\fui,\psi\in{\cD}_R$ define $d(\fui,\psi)=E(\psi^{-1}\fui)$.

If $H,K\in\cH_R$ generate $\fui,\psi$ respectively, then
$L(x,t)=(H-K)(\psi_t(x),t)$ belongs to
${\cH}_0$ and generates $\psi^{-1}\fui$. 
Reciprocally, suppose that $K\in\cH_R$ generates $\psi$. 
If $L\in{\cH}_0$ generates $\psi^{-1}\fui$ then
$H(x,t)=K(x,t)+L(\psi_t^{-1}(x),t)$ belongs to $\cH_R$
and generates $\fui$.

Thus, for any $\psi,\fui\in\cD_R$ and $K\in\cH_R$ generating $\psi$ we
have
\begin{multline*}
\left\{\lV(H-J)_t\rV:H,J\in\cH_R,\fui=\fui^H_1,\psi=\fui^J_1\right\}\\ 
\subset \left\{\lV L_t\rV:L\in{\cH}_0,\psi^{-1}\fui=\fui^L_1\right\}\\
\subset  \left\{\lV (H-K)_t\rV dt:H\in\cH_R,\fui=\fui^H_1\right\}.
\end{multline*}

Therefore
\begin{eqnarray*}
   d(\fui,\psi)  &  = & \inf\left
\{\int_0^1\lV (H-J)_t\rV dt:H,J\in\cH_R:\fui=\fui^H_1,
\psi=\fui^J_1\right\}\\
& = &\inf\left\{\int_0^1\lV (H-K)_t\rV dt:H\in\cH_R:\fui=\fui^H_1\right\}.
\end{eqnarray*}

\begin{corollary}\label{triangle}
Let $\phi\in{\cD}_0$ and define 
\[E^+(\phi)=\inf\{E^+(L):L\in{\cH}_0,\, \fui_1^L=\phi\}\]
\[E^-(\phi)=\sup\{E^-(L):L\in{\cH}_0,\, \fui_1^L=\phi\}.\]

Let $\psi\in\cD_R$, then
\[E^-(\phi)\le \ga_\pm(\psi\circ\phi)-\ga_\pm(\psi)\le E^+(\phi).\]
\end{corollary}

{\it Proof}
Let $L\in{\cH}_0$, $K\in{\cH}_R$ with $\fui_1^L=\phi$ and
$\fui_1^K=\psi$. As above $H(x,t)=K(x,t)+L(\psi_t^{-1}(x),t)$ belongs to 
$\cH_R$ and generates $\psi\circ\phi$. By item 2 in Theorem \ref{critical}
\[E^-(L)\le \ga_\pm(H)-\ga_\pm(K)\le E^+(L),\]
from which Corollary \ref{triangle} follows.
\qed

\begin{definition}
A function $H\in{\cH}_0$ is called quasi-autonomous if
there are $x_-,x_+\in U_1$ such that
\[H(x_-,t)=\min H_t,\quad H(x_+,t)=\max H_t\]
for all $t\in{\mS}^1$.
\end{definition}

Recall that $G_r(q,p)=(q,\exp_q(rp))$.
Let $\cG$ be the set of functions $S:G_r(U_{1+\de})\to{\R}$ such that 
\begin{itemize}
\item $\dfrac{\partial^2 S}{\partial q\partial Q}$ is negative definite.
\item $S(q,Q)=\dfrac{d(q,Q)^2}{2r}$ for $d(q,Q)\ge r$.
\end{itemize}

\begin{proposition}\label{ham-jac}
Let $S_t,t\in[a,b]$ be a path in $\cG$. Then the
corresponding path $\fui_t$ in $\cT$ is generated by a Hamiltonian 
$H\in{\cH}_0$ satisfying the Hamilton-Jacobi equation
\begin{equation}
  \label{eq:ham-jac}
  \frac{\partial S_t}{\partial t}(q,Q)+
H(Q,\frac{\partial S_t}{\partial Q}(q,Q),t)=0.
\end{equation}
Moreover, $H$ is quasi-autonomous if and only if
$\dfrac{\partial S_t}{\partial t}$ is. 
\end{proposition}

The proof of equation \ref{eq:ham-jac} is standard and the proof of
the last statement is the same as the given by Bialy and Polterovich
\cite{B-P}.

Consider the length function ${\cL}$ on the fundamental group of ${\cD}_0$
\[{\cL}([\fui_t])= \inf_{\psi_t\in[\fui_t]}\text{length}(\psi_t).\]

We claim that the image of ${\cL}$ is $\{0\}$. In fact, given any loop  
$\fui_t=(Q_t,P_t)$  in ${\cD}_0$ with generating Hamiltonian 
$H\in{\cH}_0$, we can define the loop $\psi:{\mS}^1\to{\cD}_0$ 
by $\psi_t(q,p)=(Q_t(q,p/\ep),\ep P_t(q,p/\ep))$
with generating Hamiltonian $F(q,p,s)=\ep H(q,p/\ep,s)$.
Therefore 
\[\text{length}(\psi)= \int_0^1\lV F_t\rV dt = \ep \int_0^1\lV H_t\rV dt.\]
Thus, for any class $[\fui]\in\pi_1({\cD}_0)$ we have
$ {\cL}([\fui])= 0$.

\begin{proposition}\label{quasi}
There is a $C^2$ neighborhood $\cU$ of zero in ${\cH}_0$ such that 
\[\lV \fui_1^H \rV = \int_0^1\lV H_t\rV dt \]
for any quasi-autonomous $H\in \cU$ .
\end{proposition}
{\it Proof}.
By Remark 3.3 and item (ii) of Theorem 1.3 in \cite{L-M} II, 
there is a $C^2$ neighborhood $\cU$ of zero in ${\cH}_0$ such that 
for any  quasi-autonomous $H\in \cU$ the path $\fui^H_{t\in[0,1]}$ 
is length-minimizing amongst all paths homotopic (with fixed end points) to 
$\fui^H_{t\in[0,1]}$.

Suppose that $\phi$ is other path in ${\cD}_0$ with the same 
end points and length$(\phi)<$ length$(\fui^H)$. 
Choose a loop $\psi\in[-\phi*\fui^H]$ such that
\[\text{length}(\psi)< \text{length}(\fui^H)-\text{length}(\phi).\]
the path $\phi*\psi$ is homotopic to $\fui^H$ and shorter: a contradiction.

\qed

{\it Proof of Proposition \ref{flat}}.
For $\fui_0,\fui_1\in\cT$, let $S_0,S_1\in{\cG}$ 
be their generating functions. For $t\in[0,1]$,
$S_t=(1-t)S_0+tS_1$ defines a map $\fui_t\in{\cT}$. 
Consider the path $\fui_t\fui_0^{-1}$ in ${\cD}_0$ and 
its generating Hamiltonian $H$. 
As in \cite{B-P}, one easily shows that 
$\dfrac{\partial S_t}{\partial t}=S_1-S_0$ is quasi-autonomous. 
By Proposition \ref{ham-jac}, $H$ is quasi-autonomous and
\[\lV H_t\rV=\lV S_0-S_1\rV.\]

By Proposition \ref{quasi}, there is a $C^1$ neighborhood $\mathcal O$
of $\fui$ such that if $\fui_0,\fui_1\in\mathcal O$,
the length of the path $\fui_t\fui_0^{-1}$ equals
$d(\fui_0,\fui_1)$.

\section{Appendix}

    Let $\te$ be a closed 1-form in $M$ such that $\te\in[\omega]$.
    Let $\Ga:[0,n]\to T^*M$  be a holonomic curve meaning that for
    $\ga=\pi\circ\Ga$ one has $\dot{\ga}(t)=\partial H/\partial p(\Ga(t),t)$.
    For each $(q,p,t)\in T^*M\times{\mS}^1$ there is $\xi\in T^*_qM$ such that
    \[H(q,\te(q),t)=H(q,p,t)+(\te(q)-p)\frac{\partial H}{\partial p}(q,p,t)
    +\frac 12\, p\,\frac{\partial^2 H}{\partial p^2}(q,\xi,t)\,p,\]

    so

    \[(p-\te(q))\frac{\partial H}{\partial p}(q,p,t)-H(q,p,t)=-H(q,\te(q),t)
   +\frac 12\, p\,\frac{\partial^2 H}{\partial  p^2}(q,\xi,t)\,p.\]

    Let $L:TM\times{\mS}^1\to\R$ be the Legendre transform of $H$.
    Since $H$ is convex, for each $n\in{\Z}^+$ we have
    \begin{multline*}
      \int_0^n (L(\ga,\dot{\ga},t)-\te(\ga)\dot{\ga})dt=
      \int_\Ga\lam-\pi^*\te-Hdt\\
      \ge-\int_0^n H(\ga(t),\te(\ga(t)),t)dt
      \ge -\int_0^n\max_{q\in M}H(q,\te(q),t)dt\\
      =- n\int_0^1\max_{q\in M}H(q,\te(q),t)dt.  
    \end{multline*}

    As the left hand side does not depend on the representant of  the
    class $[\omega]$, we have
    \[\frac 1n \int_0^n (L(\ga,\dot{\ga},t)-\omega(\ga)\dot{\ga})dt\ge 
    \sup_{\te\in[\omega]}-\int_0^1\max_{q\in M}H(q,,\te(q),t)dt=
    -c(H,[\omega]).\]

    Since this holds for any holonomic curve, it follows from Proposition
    1 in \cite{M} that
    \[-\alpha([\omega])
    =\min\left\{\int (L-\omega)\,d\mu:\;\mu\in{\M}(L)\right\}
    \ge -c(H,[\omega]).\]

    So 
    \[\alpha([\omega])\le c(H,[\omega]).\]
    \qed

    {\it Proof of Corollary \ref{equal}}.
    Since $H$ is convex, for each $t$ we have that
    \[H_t|_L=\min_{x\in U_1}H(x,t),\quad 0=\max_{x\in U_1}H(x,t)\]
    Therefore 
    \[E(\fui^H_1)\le
    \int_0^1 \lV H_t\rV dt=-\int_0^1 H_t|_L dt\le- C_H\le \beta (0)
    \le E(\fui^H_1).\]
    \qed

\end{document}